\documentclass[12pt]{article}
\usepackage{amsfonts}
\usepackage{mathrsfs}
\usepackage{amsmath,amssymb}
\numberwithin{equation}{section}

\def\c{\frak{c}}
\def\D{\Delta}
\def\fc{{\textbf{\textit c}}}

\def\LLsl{{\mathcal{L}^s_{\lambda}}}
\def\WWsl{{\mathcal{W}^s_{\lambda}}}
\def\VVsl{{\mathcal{V}^s_{\lambda}}}

\def\QQsl{{\mathcal{Q}^s_{\lambda}}}
\def\vpsl{\varphi_{\lambda}^s}

\def\LCsl{{{\mathcal{L}_\frak{C}}^s_{\lambda}}}
\def\WCsl{{{\mathcal{W}_\frak{C}}^s_{\lambda}}}
\def\Dsl{{D^s_{\lambda}}}
\def\Dslzr{{(D^s_{\lambda})_0}}
\def\Dslalp{{(D^s_{\lambda})_\alpha}}

\def\HHsl{\mathcal{H}^s_{\lambda}}
\def\JJsl{\mathcal{J}^s_{\lambda}}

\def\UU{\frak{U}}

\def\H1LL{H^1({\mathcal{L}^s_{\lambda,\mu}},{\mathcal{L}^s_{\lambda,\mu}})}

\def\Der{{\rm Der}}
\def\Inn{{\rm Inn}}
\def\Ker{{\rm Ker}}
\def\Im{{\rm Im}}

\def\sp{{\rm{span}}}

\def\cl{\centerline}
\def\ni{\noindent}
\def\vs{\vspace*}
\def\rar{\rightarrow}
\def\lrar{\longrightarrow}
\def\SM{\!\setminus\!}

\def\Z{\mathbb{Z}}
\def\C{\mathbb{C}}

\def\QED{\hfill$\Box$\par}

\newtheorem{clai}{Claim}
\newtheorem{Subclai}{Subclaim}
\newtheorem{theo}{Theorem}[section]

\newtheorem{lemm}[theo]{Lemma}
\newtheorem{prop}[theo]{Proposition}

\newtheorem{defi}[theo]{Definition}

\begin{document}

\baselineskip 18pt

\cl{{\bf Lie bialgebra structures on the deformative Schr\"{o}dinger-Virasoro algebras}
\footnote{Supported by the NSF grant 11101056 of China}}

\vs{6pt}

\cl{Huanxia Fa$^{1,2)}$}

\cl{\small $^{1)}$Wu Wen-Tsun Key Laboratory of Mathematics and School of Mathematical Sciences,}\vs{-3pt}
\cl{\small University of Science and Technology of China, Hefei 230026, China}

\cl{\small $^{2)}$School of Mathematics and Statistics, Changshu Institute
of Technology, Changshu 215500, China}

\cl{\small E-mail: sd\_huanxia@163.com}

\vs{6pt}

{\small \parskip .005 truein
\noindent{\small{\bf Abstract.}\ \ In this paper we investigate Lie bialgebra structures on the deformative Schr\"{o}dinger-Virasoro algebras mainly using the techniques introduced recently by Liu, Pei and Zhu, which indicate that all cases considered in this paper except one behave different from their centerless ones.\vs{6pt}

}


\vs{12pt}

\cl{\bf\S1. \
Introduction}\setcounter{section}{1}\setcounter{theo}{0}\setcounter{equation}{0}

Both original and twisted Schr\"{o}dinger-Virasoro Lie algebras were introduced by \cite{H1} in the context of non-equilibrium statistical physics and further investigated in \cite{H2,HU}, whose deformations were introduced in \cite{RU}. All of them are closely related to the Schr\"{o}dinger Lie algebra and the Virasoro Lie algebra. The {\it deformative Schr\"{o}dinger-Virasoro Lie algebras} $\LLsl\,(\lambda\in\C)$ considered in this paper, possess the $\C$-basis $\{L_n,\,M_n,\,Y_{s+n},\,\c\,|\,n\in\Z,\,s=0\ {\rm or}\ \frac{1}{2}\}$ with the following non-vanishing Lie brackets:
\begin{eqnarray}\begin{array}{lll}\label{LB1}
&[\,L_n,L_m]\!\!\!&=(m-n)L_{m+n}+\frac{m^3-m}{12}\delta_{m+n,0}\c,\ \ \,\,
[L_n,M_m\,\,]=(m-\lambda n)M_{m+n},\\[6pt]
&[L_n,Y_{s+m}]\!\!&=(s+m-\frac{\lambda+1}{2}n)Y_{s+m+n},\ \,\,\,\,\ \ \,\, [\,Y_{s+n},Y_{s+m}\,]=(m-n)M_{2s+m+n}.\end{array}
\end{eqnarray}
It is easy to see that $\LLsl$ contain the Virasoro algebra $\frak{v}=\sp_\C\{L_n,\,\c\,|\,n\in\Z\}$ and Witt algebra $\frak{w}=\sp_\C\{L_n\,|\,n\in\Z\}$ as their subalgebras and $\frak{h}^s=\{M_n,\,Y_{s+n},\,\c\,|\,n\in\Z\}$ as their ideals.

Due to their mathematical and physical backgrounds, interests and importance, a series of papers on Lie algebras of these types  appeared. The Lie bialgebra structures, modules of intermediate series, non-trivial vertex algebra representations and Whittaker modules over the original Schr\"{o}dinger-Virasoro Lie algebra were respectively investigated in \cite{HLS,LS-JMP,U,ZTL}. The derivations, central extensions and automorphism group of the extended Schr\"{o}dinger-Virasoro Lie algebra were determined in \cite{GJP-AC}. The Verma modules and derivations over generalized Schr\"{o}dinger-Virasoro algebras were respectively investigated in \cite{TZ-JA09} and \cite{WLXin}. The second cohomology group of both original and twisted deformative Schr\"{o}dinger-Virasoro algebras were completely determined by \cite{LJ,LSZ,RU}. Automorphisms and derivations of twisted
and original deformative Schr\"{o}dinger-Virasoro Lie algebras were respectively determined in \cite{WLXu} and \cite{JW-arXiv12}.

The notion of Lie bialgebras was introduced in 1983 by \cite{D1} to
find solutions of the Yang-Baxter equation. Many types of Lie bialgebras were considered over different algebras. Initially, the Witt and Virasoro Lie bialgebras were  investigated in \cite{T} and classified in \cite{NT}. Lie bialgebras of generalized Witt types, generalized Virasoro-like type, generalized Weyl type, Hamiltonian type and Block type were considered respectively in \cite{SS,WSS06,Y-XS08,XSS,LSX}. Lie superbialgebra structures on the Ramond $N\!=\!2$ super Virasoro algebra were considered in \cite{Y-HS09}. The Lie bialgebra structures on the original case ${\mathcal{L}^\frac12_{0}}$ (with $\c=0$) and a twisted case ${\mathcal{L}^0_{-3}}$ (with $\c=0$) were respectively investigated by \cite{HLS} and \cite{FLL}. In this paper we shall investigate Lie bialgebra structures on $\LLsl$ mainly using the techniques introduced recently in \cite{LPZ-JA12}.\vs{12pt}

\cl{\bf\S2. \
Preliminaries and main results}\setcounter{section}{2}\setcounter{theo}{0}\setcounter{equation}{0}\vs{8pt}

For convenience, some basic concepts about Lie bialgebras are presented firstly. Respectively denote $\xi$ and $\tau$ the {\it cyclic map} of $L\otimes
L\otimes L$ and the {\it twist map} of $L\otimes L$ for any given vector space $L$, which imply $ \xi(x_{1}\otimes x_{2}\otimes x_{3})=x_{2}\otimes x_{3}\otimes x_{1}$ and $\tau(x_{1}\otimes x_{2})=x_{2}\otimes x_{1}$, $\forall$ $x_1,x_2,x_3\in L$. For a vector space $L$ and two bilinear maps $\delta:L\otimes L\rar L$ and $\D: L\to L\otimes L$, the pair $(L,\delta)$ becomes a {\it Lie algebra} if the following two conditions hold:
\begin{eqnarray*}
&&\Ker(1-\tau)\subset\Ker\,\delta,\ \ \ \delta\cdot(1\otimes\delta)\cdot(1+\xi+\xi^{2})=0,
\end{eqnarray*}
and the pair $(L,\Delta)$ becomes a {\it Lie coalgebra} if the following two conditions hold:
\begin{eqnarray*}
&&\Im\,\Delta\subset\Im(1-\tau),\ \ \ (1+\xi+\xi^{2})\cdot(1
\otimes \Delta)\cdot\D =0.
\end{eqnarray*}
The triple $( L,\delta,\Delta)$ will become a {\it Lie bialgebra} if the Lie algebra $(L,\delta)$ and Lie coalgebra $(L,\Delta)$ satisfy the following compatible condition:
\begin{eqnarray}
\label{compc}&&\Delta\delta (x\otimes y)=x\cdot\Delta y-y\cdot \Delta x,\ \
\forall\,\,x,y\in L,
\end{eqnarray}
where the symbol ``$\cdot$'' is defined to be the {\it diagonal adjoint action}:
\begin{eqnarray*}
&&x\cdot(\mbox{$\sum\limits_{i}$}{a_{i}\otimes b_{i}})=
\mbox{$\sum\limits_{i}$}({[x,a_{i}]\otimes
b_{i}+a_{i}\otimes[x,b_{i}]}).
\end{eqnarray*}

Denote $\UU$ the universal enveloping algebra of $L$ and ${\bf 1}$ the
identity element of $\UU$. For any $r =\sum_{i} {a_{i} \otimes
b_{i}}\in L\otimes L$, we denote
\begin{eqnarray*}
r^{12}=\mbox{$\sum\limits_{i}$}{a_{i} \otimes b_{i} \otimes{\bf 1}},\ \ r^{13}=
\mbox{$\sum\limits_{i}$}{a_{i} \otimes{\bf 1}\otimes b_{i}},\ \ r^{23}=\mbox{$\sum\limits_{i}$}{\bf 1}
\otimes a_{i} \otimes b_{i}.
\end{eqnarray*}
Define $\fc(r)$ to be elements of $\UU
\otimes \UU \otimes \UU$ by
\begin{eqnarray*}
\fc(r)\!\!\!&=&\!\!\![r^{12},r^{13}]+[r^{12},r^{23}]+[r^{13},r^{23}]\\
\!\!\!&=&\!\!\!\mbox{$\sum\limits_{i,j}$}[a_i,a_j]\otimes b_i\otimes
b_j+\mbox{$\sum\limits_{i,j}$}a_i\otimes [b_i,a_j]\otimes b_j+
\mbox{$\sum\limits_{i,j}$}a_i\otimes a_j\otimes [b_i,b_j].
\end{eqnarray*}

\begin{defi}\label{def1}

(1) A {\it coboundary Lie bialgebra} is a $4$-tuple $( L,
\delta, \D,r),$ where $( L,\delta,\D)$ is a Lie bialgebra and $r \in
\Im(1-\tau)\subset L\otimes L$ such that $\D=\D_r$ is a coboundary of $r$, i.e.,
\begin{eqnarray}
\label{D-r}\D_r(x)=x\cdot r\mbox{\ \ for\ \ }x\in L.
\end{eqnarray}

(2)\ \ A coboundary Lie bialgebra $( L,\delta,\D,r)$ is called {\it
triangular} if it satisfies the following classical Yang-Baxter
Equation
\begin{eqnarray}\label{CYBE}
\fc(r)=0.
\end{eqnarray}

(3)\ \ $r\in\Im(1-\tau)\subset L\otimes L$ is said to
satisfy the modified Yang-Baxter equation if
\begin{eqnarray}\label{MYBE}
x\cdot \fc(r)=0,\ \,\forall\,\,x\in L.
\end{eqnarray}
\end{defi}

Regard $\VVsl=\LLsl\otimes\LLsl$ as an $\LLsl$-module under the adjoint
diagonal action. Denote by $\Der(\LLsl,\VVsl)$ the set of
\textit{derivations} $\Dsl:\LLsl\to\VVsl$, namely, $\Dsl$ is a linear map
satisfying
\begin{eqnarray}\label{deriv}
\Dsl([x,y])=x\cdot\Dsl(y)-y\cdot\Dsl(x),
\end{eqnarray}
and $\Inn(\LLsl,\VVsl)$ the set consisting of the derivations $v_{\rm
inn}, \, v\in\VVsl$, where $v_{\rm inn}$ is the \textit{inner
derivation} defined by $v_{\rm inn}:x\mapsto x\cdot v.$ Then $H^1(\LLsl,\VVsl)\cong\Der(\LLsl,\VVsl)/\Inn(\LLsl,\VVsl)$,
where $H^1(\LLsl,\VVsl)$ is the {\it first cohomology group} of the Lie
algebra $\LLsl$ with coefficients in the $\LLsl$-module $\VVsl$.\\

Lemmas \ref{Origderivnstns} and \ref{Twistderivnstns} can be found in \cite{JW-arXiv12}, \cite{LS-arXiv08} and \cite{WLXu}.
\begin{lemm}\label{Origderivnstns}
\begin{eqnarray*}
\begin{array}{lll}
&&H^1({\mathcal{L}^\frac12_{\lambda}},{\mathcal{L}^\frac12_{\lambda}})
\cong\left\{\begin{array}{lllll}
\C D_1\oplus\C D^{0}_2\oplus\C D^{0}_3&\mbox{\rm{if}\ \ }\lambda=0,\\[3pt]
\C D_1\oplus\C D^{-1}_2&\mbox{\rm{if}\ \ }\lambda=-1,\\[3pt]
\C D_1\oplus\C D^{-2}_2&\mbox{\rm{if}\ \ }\lambda=-2,\\[3pt]
\C D_1&\mbox{{\rm otherwise}},
\end{array}\right.
\end{array}
\end{eqnarray*}
where the corresponding non-vanishing components are given as follows:
\begin{eqnarray*}
D_1(Y_{n+\frac12})\!\!\!&=&\!\!\!Y_{n+\frac12},\ \ \ D_1(M_n)=2M_n,\ \ \ D^{0}_2(L_n)=M_n,\ \ \ D^{0}_3(L_n)=nM_n,\\[3pt]
D^{-1}_2(L_n)\!\!\!&=&\!\!\!(n^2-n)M_n,\ \ \ D^{-2}_2(L_n)=n^3M_n,
\end{eqnarray*}
for all $n\in\Z$.
\end{lemm}

\begin{lemm}\label{Twistderivnstns}
\begin{eqnarray*}
\begin{array}{lll}
&&H^1({\mathcal{L}^0_{\lambda}},{\mathcal{L}^0_{\lambda}})
\cong\left\{\begin{array}{lllll}
\C d_1\oplus\C d^{0}_2\oplus\C d^{0}_3&\mbox{\rm{if}\ \ }\lambda=0,\\[3pt]
\C d_1\oplus\C d^{-1}_2\oplus\C d^{-1}_3&\mbox{\rm{if}\ \ }\lambda=-1,\\[3pt]
\C d_1\oplus\C d^{-2}_2&\mbox{\rm{if}\ \ }\lambda=-2,\\[3pt]
\C d_1\oplus\C d^{1}_2&\mbox{\rm{if}\ \ }\lambda=1,\\[3pt]
\C d_1&\mbox{{\rm otherwise}},
\end{array}\right.
\end{array}
\end{eqnarray*}
where the corresponding non-vanishing components are given as follows:
\begin{eqnarray*}
d_1(Y_{n})\!\!\!&=&\!\!\!Y_{n},\ \ \ d_1(M_n)=2M_n,\ \ \ d^{0}_2(L_n)=M_n,\ \ \
d^{0}_3(L_n)=nM_n,\\[3pt]
d^{-1}_2(L_n)\!\!\!&=&\!\!\!n^2M_n,\ \ \
d^{-1}_3(Y_n)=nM_n,\ \ \ d^{-2}_2(L_n)=n^3M_n,\ \ \ d^{1}_2(Y_{n})=M_n,
\end{eqnarray*}
for all $n\in\Z$.
\end{lemm}

The Lie bialgebra structures on the Schr\"{o}dinger-Virasoro Lie algebras $\mathcal{L}^0_{-3}$ and $\mathcal{L}^\frac12_{0}$ have already been considered respectively in \cite{FLL} and \cite{HLS} (with $\c=0$). We just need to consider the following cases:
\begin{eqnarray}
&&\mathcal{L}^\frac12_{-2},\ \ \ \mathcal{L}^\frac12_{-1},\ \ \ \mathcal{L}^\frac12_{\lambda}\ \ \big(\,\forall\,\,\lambda\notin \mathcal{S}_\frac12\big),\label{3--12cases}\\[3pt]
&&\mathcal{L}^0_{-2},\ \ \ \mathcal{L}^0_{-1},\ \ \ \mathcal{L}^0_{0},\ \ \ \mathcal{L}^0_{1},\ \ \ \mathcal{L}^0_{\lambda}\ \ \big(\,\forall\,\, \lambda\notin\mathcal{S}_0\big),\label{5--0cases}
\end{eqnarray}
where $\mathcal{S}_\frac12\!=\!\{0,\,-1,\,-2\}$ and $\mathcal{S}_0\!=\!\{0,\,\pm1,\,-2,\,-3\}$.

Lemmas \ref{Origderivnstns} and \ref{Twistderivnstns} can be partially rewritten as follows.
\begin{lemm}\label{RewOrigderivnstns}
\begin{eqnarray*}
\begin{array}{lll}
&&H^1({\mathcal{L}^\frac12_{\lambda}},{\mathcal{L}^\frac12_{\lambda}})
\cong\left\{\begin{array}{lllll}
\frak{D}_1&\mbox{\rm{if}\ \ }\lambda=-1,\\[3pt]
\frak{D}_2&\mbox{\rm{if}\ \ }\lambda=-2,\\[3pt]
\frak{D}_3&\mbox{\rm{if}\ \ }\lambda\notin\mathcal{S}_\frac12,
\end{array}\right.
\end{array}
\end{eqnarray*}
where $\frak{D}_i$ consist of the following derivations $\sigma_i$ respectively:
\begin{eqnarray*}
\sigma_1(L_n)\!\!\!&=&\!\!\!\alpha_1(n^2-n)M_n,\ \ \ \sigma_1(Y_{n+\frac12})=\beta_1Y_{n+\frac12},\ \ \ \sigma_1(M_n)=2\beta_1M_n,\\[3pt]
\sigma_2(L_n)\!\!\!&=&\!\!\!\alpha_2n^3M_n,\ \ \ \sigma_2(Y_{n+\frac12})=\beta_2Y_{n+\frac12},\ \ \ \sigma_2(M_n)=2\beta_2M_n,\\[3pt]
\sigma_3(Y_{n+\frac12})\!\!\!&=&\!\!\!\beta_3Y_{n+\frac12},\ \ \ \sigma_3(M_n)=2\beta_3M_n,
\end{eqnarray*}
for all $n\in\Z$ and any $\alpha_1,\,\alpha_2,\,\beta_i\in\C$ with $i=1,\,2,\,3$.
\end{lemm}

\begin{lemm}\label{RewTwistderivnstns}
\begin{eqnarray*}
\begin{array}{lll}
&&H^1({\mathcal{L}^0_{\lambda,\mu}},{\mathcal{L}^0_{\lambda,\mu}})
\cong\left\{\begin{array}{lllll}
\frak{d}_1&\mbox{\rm{if}\ \ }\lambda=0,\\[3pt]
\frak{d}_2&\mbox{\rm{if}\ \ }\lambda=-1,\\[3pt]
\frak{d}_3&\mbox{\rm{if}\ \ }\lambda=-2,\\[3pt]
\frak{d}_4&\mbox{\rm{if}\ \ }\lambda=1,\\[3pt]
\frak{d}_5&\mbox{\rm{if}\ \ }\lambda\notin\mathcal{S}_0,
\end{array}\right.
\end{array}
\end{eqnarray*}
where $\frak{d}_i$ consist of the following derivations $\varrho_i$ respectively:
\begin{eqnarray*}
\varrho_1(L_n)\!\!\!&=&\!\!\!(\mu_1n+\nu_1)M_n,\ \ \
\varrho_1(Y_{n})=\gamma_1Y_{n},\ \ \ \varrho_1(M_n)=2\gamma_1M_n,\\[3pt]
\varrho_2(L_n)\!\!\!&=&\!\!\!\mu_2n^2M_n,\ \ \
\varrho_2(Y_{n})=\gamma_2Y_{n}+\zeta_2\,nM_n,\ \ \ \varrho_2(M_n)=2\gamma_2M_n,\\[3pt]
\varrho_3(L_n)\!\!\!&=&\!\!\!\mu_3n^3M_n,\ \ \ \varrho_3(Y_{n})=\gamma_3Y_{n},\ \ \ \varrho_3(M_n)=2\gamma_3M_n,\\[3pt]
\varrho_4(Y_{n})\!\!\!&=&\!\!\!\gamma_4Y_{n}+\zeta_4M_n,\ \ \ \varrho_4(M_n)=2\gamma_4M_n,\\[3pt]
\varrho_5(Y_{n})\!\!\!&=&\!\!\!\gamma_5Y_{n},\ \ \ \varrho_5(M_n)=2\gamma_5M_n,
\end{eqnarray*}
for all $n\in\Z$ and any $\mu_1,\,\mu_2,\,\mu_3,\,\nu_1,\,\zeta_2,\,\zeta_4,\,\gamma_i\in\C$ with $i=1,\,2,\,3,\,4,\,5$.
\end{lemm}

A derivation $\Dsl\in\Der(\LLsl,\VVsl)$ is {\it homogeneous of degree $\alpha\in\Z_s$} if $\Dsl\big((\LLsl)_p\big)\subset(\VVsl)_{\alpha+p}$ for all $p\in\Z_s$. Denote $\Der(\LLsl,\VVsl)_\alpha =\{\Dsl\in\Der(\LLsl,\VVsl)\,|\,{\rm deg\,}\Dsl=
\alpha\}$ for $\alpha\in\Z_s$.

Let $\Dsl$ be an element of
$\Der(\LLsl,\VVsl)$. For any $\alpha\in\Z_s$, define the linear map
$\Dslalp:\LLsl\rightarrow\VVsl$ as follows: For any $\mu\in(\LLsl)_q$ with
$q\in\mathbb{Z}_s$, write $\Dsl(\mu)=\sum_{p\in\mathbb{Z}_s}\mu_p$ with
$\mu_p\in(\VVsl)_p$, then we set $\Dslalp(\mu)=\mu_{q+\alpha}$.
Obviously, $\Dslalp\in \Der(\LLsl,\VVsl)_\alpha$ and we have
\begin{eqnarray}\label{summable}
\Dsl=\mbox{$\sum\limits_{\alpha\in\mathbb{Z}_s}\Dslalp$},
\end{eqnarray}
which holds in the sense that for every $u\in\LLsl$, only finitely
many $\Dslalp(u)\neq 0,$ and
$\Dsl(u)=\sum_{\alpha\in\mathbb{Z}_s}\Dslalp(u)$ (we call such a sum in
(\ref{summable}) {\it summable}).

\begin{lemm}\label{OrderTensor}
One can find the following elements of $\Der({\mathcal{L}^\frac12_{\lambda}},
{\mathcal{L}^\frac12_{\lambda}}\otimes{\mathcal{L}^\frac12_{\lambda}})$: $\sigma^\natural_1$ for the case $\lambda=-1$,
$\sigma^\natural_2$ for the case $\lambda=-2$,
$\sigma^\natural_3$ for the case $\lambda\notin\mathcal{S}_\frac12$,
which are defined by the following relations respectively:
\begin{eqnarray*}
\sigma^\natural_1(L_n)\!\!\!&=&\!\!\!(n^2-n)(\alpha_1z_1\otimes M_n+\alpha^\dag_1M_n\otimes z^\dag_1),\\[3pt]
\sigma^\natural_1(Y_{n+\frac12})\!\!\!&=&\!\!\!\beta_1w_1\otimes Y_{n+\frac12}+\beta^\dag_1Y_{n+\frac12}\otimes w^\dag_1,\ \ \ \sigma^\natural_1(M_n)=2(\beta_1w_1\otimes M_n+\beta^\dag_1M_n\otimes w^\dag_1),\\[3pt]
\sigma^\natural_2(L_n)\!\!\!&=&\!\!\!n^3(\alpha_2z_2\otimes M_n+\alpha^\dag_2M_n\otimes z^\dag_2),\ \ \ \sigma^\natural_2(Y_{n+\frac12})=\beta_2w_2\otimes Y_{n+\frac12}+\beta^\dag_2Y_{n+\frac12}\otimes w^\dag_2,\\[3pt] \sigma^\natural_2(M_n)\!\!\!&=&\!\!\!2(\beta_2w_2\otimes M_n+\beta^\dag_2M_n\otimes w^\dag_2),\ \ \ \sigma^\natural_3(Y_{n+\frac12})=\beta_3w_3\otimes Y_{n+\frac12}+\beta^\dag_3 Y_{n+\frac12}\otimes w^\dag_3,\\[3pt]
\sigma^\natural_3(M_n)\!\!\!&=&\!\!\!2(\beta_3w_3\otimes M_n+\beta^\dag_3 M_n\otimes w^\dag_3),\ \ \ \sigma^\natural_3(L_n)=\sigma^\natural_1(\c)=\sigma^\natural_2(\c)=\sigma^\natural_3(\c)=0,
\end{eqnarray*}
for all $n\in\Z$ and some $\alpha_1,\,\alpha_2,\,\alpha^\dag_1,\,\alpha^\dag_2,\,\beta_i\in\C$ and $w_i,\,z_j,\,w^\dag_i,\,z^\dag_j\in\frak{C}$ with $i=1,\,2,\,3$, $j=1,\,2$.
\end{lemm}

\begin{lemm}\label{TwderTensor}
One can find the following elements of $\Der({\mathcal{L}^0_{\lambda}},
{\mathcal{L}^0_{\lambda}}\otimes{\mathcal{L}^0_{\lambda}})$: $\varrho^\natural_1$ for the case $\lambda=0$, $\varrho^\natural_2$ for the case $\lambda=-1$,  $\varrho^\natural_3$ for the case $\lambda=-2$, $\varrho^\natural_4$ for the case $\lambda=1$, $\varrho^\natural_5$ for the case $\lambda\notin\mathcal{S}_0$,
defined by the following relations respectively:
\begin{eqnarray*}
\varrho^\natural_1(L_n)\!\!\!&=&\!\!\!(\mu_1n+\nu_1)z_1\otimes M_n+(\mu^\dag_1n+\nu^\dag_1)M_n\otimes z^\dag_1,\ \ \
\varrho^\natural_1(Y_{n})=\gamma_1w_1\otimes Y_{n}+\gamma^\dag_1Y_n\otimes w^\dag_1,\\[3pt] \varrho^\natural_1(M_n)\!\!\!&=&\!\!\!2(\gamma_1w_1\otimes M_n+\gamma^\dag_1M_n\otimes w^\dag_1),\ \ \
\varrho^\natural_2(L_n)=n^2(\mu_2z_2\otimes M_n+\mu^\dag_2M_n\otimes z^\dag_2),\\[3pt]
\varrho^\natural_2(Y_{n})\!\!\!&=&\!\!\!\gamma_2w_2\otimes Y_{n}+\gamma^\dag_2Y_n\otimes w^\dag_2+n(\zeta_2v_2\otimes M_n+\zeta^\dag_2M_n\otimes v^\dag_2),\\[3pt] \varrho^\natural_2(M_n)\!\!\!&=&\!\!\!2(\gamma_2w_2\otimes M_n+\gamma^\dag_2M_n\otimes w^\dag_2),\ \ \ \varrho^\natural_3(L_n)=n^3(\mu_3z_3\otimes M_n+\mu^\dag_3M_n\otimes z^\dag_3),\\[3pt]
\varrho^\natural_3(Y_{n})\!\!\!&=&\!\!\!\gamma_3w_3\otimes Y_{n}+\gamma^\dag_3Y_n\otimes w^\dag_3,\ \ \ \varrho^\natural_3(M_n)=2(\gamma_3w_3\otimes M_n+\gamma^\dag_3M_n\otimes w^\dag_3),\\[3pt]
\varrho^\natural_4(Y_{n})\!\!\!&=&\!\!\!\gamma_4w_4\otimes Y_{n}+\gamma^\dag_4 Y_{n}\otimes w^\dag_4+\zeta_4v_4\otimes M_n+\zeta^\dag_4M_n\otimes v^\dag_4,\\[3pt] \varrho^\natural_4(M_n)\!\!\!&=&\!\!\!2(\gamma_4w_4\otimes M_n+\gamma^\dag_4 M_{n}\otimes w^\dag_4),\ \ \
\varrho^\natural_5(Y_{n})=\gamma_5w_5\otimes Y_{n}+\gamma^\dag_5 Y_{n}\otimes w^\dag_5,\\[3pt]
\varrho^\natural_5(M_n)\!\!\!&=&\!\!\!2(\gamma_5w_5\otimes M_n+\gamma^\dag_5 M_{n}\otimes w^\dag_5),\ \ \ \varrho^\natural_4(L_n)=\varrho^\natural_5(L_n)=\varrho^\natural_k(\c)=0,
\end{eqnarray*}
for all $n\in\Z$, some $\mu_i,\,\mu^\dag_i,\,\nu_1,\,\nu^\dag_1,\,\zeta_2,\,\zeta_4,\,\zeta^\dag_2,\,\zeta^\dag_4,
\,\gamma_j,\,\gamma^\dag_j\in\C$ with $i=1,\,2,\,3$, $j=1,\,2,\,3,\,4,\,5$ and $z_1,\,z^\dag_1,\,w_1,\,w^\dag_1\in\frak{C}_0$, $v_2,\,v^\dag_2,\,v_4,\,v^\dag_4,\,z_k,\,z^\dag_k,\,w_k,\,w^\dag_k\in\frak{C}$ with $k=2,\,3,\,4,\,5$.
\end{lemm}

Respectively, denote the vector spaces spanned by $\sigma^\natural_i$ and $\varrho^\natural_j$ as $\frak{D}^\natural_i$ and $\frak{d}^\natural_j$ for $i=1,\,2,\,3$ and $j=1,\,2,\,3,\,4,\,5$.

The main results of this paper can be formulated as follows.
\begin{theo}\label{maintheo}
(i)\ \ Every Lie bialgebra $(\LLsl,[\cdot,\cdot],\D)$ is triangular coboundary for the Lie algebras $\LLsl$ given in \eqref{3--12cases} and \eqref{5--0cases} with $\c\!=\!0$ and $(\lambda,s)\neq(0,0)$.

(ii)\ \ Every Lie bialgebra $({\mathcal{L}^0_{0}},[\cdot,\cdot],\D)$ is not triangular coboundary with $\c\!=\!0$ and $\varrho^\natural_1\neq0$.

(iii)\ \ No Lie bialgebra $(\LLsl,[\cdot,\cdot],\D)$ with $\c\!\neq\!0$ is triangular coboundary if the derivations given in Lemmas \ref{OrderTensor} and \ref{TwderTensor} are not equal to zero.
\end{theo}

\vs{18pt}\cl{\bf\S3\ \ Proof of the main result}\setcounter{section}{3}\setcounter{theo}{0}\setcounter{equation}{0}\vs{10pt}

Throughout the paper we denote by $\Z^*$ the set of all nonnegative
integers and $\C^*$ the set of all nonnegative complex numbers.

The following lemma can be found in Lemma 2.2 of \cite{LPZ-JA12}.
\begin{lemm}\label{Zcenterlamm}
Suppose that $\frak{g}=\oplus_{n\in\Z}\frak{g}_n$ is a $\Z$-graded Lie algebra with a finite-dimensional center $\frak{C_g}$, and $\frak{g}_0$ is generated by $\{\frak{g}_n,\,n\neq0\}$. Then
\begin{eqnarray*}
H^1(\frak{g},\frak{C_g}\otimes\frak{g}+\frak{g}\otimes\frak{C_g})_0
=\frak{C_g}\otimes H^1(\frak{g},\frak{g})_0+H^1(\frak{g},\frak{g})_0\otimes\frak{C_g}.
\end{eqnarray*}
\end{lemm}
Denote $\Z_s=\Z\cup\{s+\Z\}$ and $\Z^*_s=\Z_s\SM\{0\}$, i.e., $\Z_0=\Z$, $\Z^*_0=\Z^*$, $\Z_\frac12=\frac12\Z$ and $\Z^*_\frac12=\frac12\Z^*$. Then the above lemma can be generalized immediately as follows.
\begin{lemm}\label{12Zcenterlamm}
Suppose that $\frak{g}=\oplus_{n\in\Z_s}\frak{g}_n$ is a $\Z_s$-graded Lie algebra with a finite-dimensional center $\frak{C_g}$, and $\frak{g}_0$ is generated by $\{\frak{g}_n,\,n\neq0\}$. Then
\begin{eqnarray*}
H^1(\frak{g},\frak{C_g}\otimes\frak{g}+\frak{g}\otimes\frak{C_g})_0
=\frak{C_g}\otimes H^1(\frak{g},\frak{g})_0+H^1(\frak{g},\frak{g})_0\otimes\frak{C_g}.
\end{eqnarray*}
\end{lemm}

The following lemma is one of the main results given in \cite{NT}.
\begin{lemm}\label{WittVirLiebialg}
Every Lie bialgebra on the Witt algebra $\frak{w}$ and Virasoro algebra $\frak{v}$ is triangular coboundary and $H^1(\frak{w},\frak{w}\otimes\frak{w})=H^1(\frak{v},\frak{v}\otimes\frak{v})=0$.
\end{lemm}

Denote $\LLsl^{\otimes n}$ the tensor product of $n$ copies of
$\LLsl$ and regard it as an $\LLsl$-module under the adjoint diagonal action of
$\LLsl$. The first item of the following lemma can be obtained by using the similar arguments as those given in the known references and the left two can be found in the references (e.g. \cite{D1,D2,NT,WSS06}). For convenience, we introduce the following notations: $\frak{C}=\C\c$,  $\frak{C}_0=\sp_\C\{\c,\,M_0\}$ and $\frak{C}^s_{\lambda}$, where $\frak{C}^s_{\lambda}=\frak{C}_0$ for the case $(s,\lambda)=(0,0)$ and $\frak{C}^s_{\lambda}=\frak{C}$ for all the other cases referred in \eqref{3--12cases} and \eqref{5--0cases}.
\begin{lemm}\label{center-Cohom}
\begin{itemize}
\item[\rm(i)]
If $x\cdot r=0$ for some $r\in\LLsl^{\otimes n}$ and all $x\in\LLsl$, then $r\in {\frak{C}^s_{\lambda}}^{\otimes n}$.
\item[\rm(ii)] The $r$ satisfies (\ref{CYBE}) if and
only if it satisfies (\ref{MYBE}).
\item[\rm(iii)]
Let $L$ be a Lie
algebra and $r\in\Im(1-\tau)\subset L\otimes L,$\  then
\begin{eqnarray*}
(1+\xi+\xi^{2})\cdot(1\otimes\D_r)\cdot\D_r(x)=x\cdot
\fc(r),\ \ \forall\,\,x\in L,
\end{eqnarray*}
and the triple $(L,[\cdot,\cdot], \D_r)$ is a Lie bialgebra if and
only if $r$ satisfies (\ref{CYBE}). \end{itemize}
\end{lemm}

\begin{prop}\label{propCohLlsVls}
\begin{eqnarray*}
\begin{array}{lll}
&&H^1(\LLsl,\VVsl)=\Der(\LLsl,\VVsl)/\mathrm{Inn}(\LLsl,\VVsl)\cong
\left\{\begin{array}{lllll}
\frak{D}^\natural_1&\mbox{\rm{if}\ \ }s=\frac12,\ \,\lambda=-1,\\[3pt]
\frak{D}^\natural_2&\mbox{\rm{if}\ \ }s=\frac12,\ \,\lambda=-2,\\[3pt]
\frak{D}^\natural_3&\mbox{\rm{if}\ \ }s=\frac12,\ \,\lambda\notin\mathcal{S}_\frac12,\\[3pt]
\frak{d}^\natural_1&\mbox{\rm{if}\ \ }s=0,\ \,\lambda=0,\\[3pt]
\frak{d}^\natural_2&\mbox{\rm{if}\ \ }s=0,\ \,\lambda=-1,\\[3pt]
\frak{d}^\natural_3&\mbox{\rm{if}\ \ }s=0,\ \,\lambda=-2,\\[3pt]
\frak{d}^\natural_4&\mbox{\rm{if}\ \ }s=0,\ \,\lambda=1,\\[3pt]
\frak{d}^\natural_5&\mbox{\rm{if}\ \ }s=0,\ \,\lambda\notin\mathcal{S}_0.
\end{array}\right.
\end{array}
\end{eqnarray*}
\end{prop}
{\bf Proof of Proposition \ref{propCohLlsVls}}\ \ \,This proposition follows from a series of claims.

Denote $\HHsl=\LLsl\otimes\frak{h}^s+\frak{h}^s\otimes\LLsl$. Then $\HHsl$ is an $\LLsl$-submodule of $\VVsl$, since $\frak{h}^s$ is an ideal of $\LLsl$ and denote the quotient $\LLsl$-module $\VVsl/\HHsl$ as $\QQsl$, on which $\frak{h}^s$ acts trivially. The exact sequence $0\rar\HHsl\rar\VVsl\rar\VVsl/\HHsl\rar0$ induces the following long exact sequence
\begin{eqnarray*}
\lrar H^0(\LLsl,\QQsl)\lrar H^1(\LLsl,\HHsl)\lrar H^1(\LLsl,\VVsl)\lrar H^1(\LLsl,\QQsl)\lrar
\end{eqnarray*}
of $\Z_s$-graded vector spaces, where all coefficients of the tensor products are in $\C$. It is easy to see that $H^0(\LLsl,\QQsl)=\QQsl^{\LLsl}=\{x\in\QQsl\,|\LLsl\cdot x=0\,\}=0$. Then $H^1(\LLsl,\HHsl)\cong H^1(\LLsl,\VVsl)$ if we can prove $H^1(\LLsl,\QQsl)=0$.

Denote $\LCsl=\LLsl\otimes\frak{C}^s_{\lambda}+\frak{C}^s_{\lambda}\otimes\LLsl$. Then $\LCsl$ is an $\LLsl$-submodule of $\HHsl$. The exact sequence $0\rar\LCsl\rar\HHsl\rar\HHsl/\LCsl\rar0$ induces the following long exact sequence
\begin{eqnarray*}
\lrar H^0(\LLsl,\HHsl/\LCsl)\lrar H^1(\LLsl,\LCsl)\lrar H^1(\LLsl,\HHsl)\lrar H^1(\LLsl,\HHsl/\LCsl)\lrar.
\end{eqnarray*}
It is easy to see that $H^0(\LLsl,\HHsl/\LCsl)=(\HHsl/\LCsl)^{\LLsl}=\{x\in\HHsl/\LCsl\,|\LLsl\cdot x=0\,\}=0$. Then $H^1(\LLsl,\LCsl)\cong H^1(\LLsl,\HHsl)$ if we can prove $H^1(\LLsl,\HHsl/\LCsl)=0$.

\begin{clai}\label{Dslalpnotzero}
If $\alpha\in\Z^*_s$, then
$\Dslalp\in\Inn(\LLsl,\VVsl)$.
\end{clai}

For $\alpha\neq 0$, denote
$\gamma=\alpha^{-1}\Dslalp(L_0)\in(\VVsl)_{\alpha}$. For any
$L_n,\,M_n\in(\LLsl)_{n}$ and $Y_{s+n}\in(\LLsl)_{s+n}$, applying $\Dslalp$ to $[L_0,L_n]=nL_n$, $[L_0,M_n]=nM_{n}$ and $[L_0,Y_{s+n}]=(s+n)Y_{s+n}$, we obtain \big(recalling $\Dslalp(L_n),\,\Dslalp(M_n)\in(\VVsl)_{n+\alpha}$ and $\Dslalp(Y_{s+n})\in(\VVsl)_{s+n+\alpha}$, i.e., $L_0\cdot \Dslalp(L_n)=(\alpha+n)\Dslalp(L_n)$, $L_0\cdot\Dslalp(M_n)=(\alpha+n)\Dslalp(M_n)$ and $L_0\cdot\Dslalp(Y_{s+n})=(\alpha+s+n)\Dslalp(Y_{s+n})$\big)
\begin{eqnarray*}
&&(\alpha+n)\Dslalp(L_n)-L_n\cdot\Dslalp(L_0)=n\Dslalp(L_n),\\
&&(\alpha+n)\Dslalp(M_n)-M_n\cdot\Dslalp(L_0)=n\Dslalp(M_n),\\
&&(\alpha+s+n)\Dslalp(Y_{s+n})-Y_{s+n}\cdot\Dslalp(L_0)=(s+n)\Dslalp(Y_{s+n}).
\end{eqnarray*}
 Then Claim \ref{Dslalpnotzero} follows.

\begin{clai}\label{DslzrL0}
$\Dslzr(L_0)\equiv0\,\,({\rm mod}\,\,\frak{C}^s_{\lambda}\otimes\frak{C}^s_{\lambda})$.
\end{clai}

For any
$L_n,\,M_n\in(\LLsl)_{n}$ and $Y_{s+n}\in(\LLsl)_{s+n}$, applying $\Dslzr$ to $[L_0,L_n]=nL_n$, $[L_0,M_n]=nM_{n}$ and $[L_0,Y_{s+n}]=(s+n)Y_{s+n}$, we obtain
\begin{eqnarray*}
&&L_n\cdot\Dslzr(L_0)=M_n\cdot\Dslzr(L_0)=Y_{s+n}\cdot\Dslzr(L_0)\equiv0\,\,({\rm mod}\,\,\frak{C}^s_{\lambda}\otimes\frak{C}^s_{\lambda}),
\end{eqnarray*}
which forces
$\Dslzr(L_0)\equiv0\,\,({\rm mod}\,\,\frak{C}^s_{\lambda}\otimes\frak{C}^s_{\lambda})$ according to Lemma \ref{center-Cohom}.

\begin{clai}\label{finsumclaim}
For any $\Dslalp\in\Der(\LLsl,\VVsl)$, \eqref{summable} is a finite sum.
\end{clai}

For any $\alpha\in\Z^*_s$, one can suppose $\Dslalp=(v_\alpha)_{\rm inn}$
for some $v_\alpha\in(\VVsl)_\alpha$. If
$\Delta_s=\{\alpha\in\Z^*_s\,|\,v_\alpha\neq0\}$ is an infinite set, then
$\Dsl(L_0)=\Dslzr(L_0)+\sum_{\alpha\in\Delta_s}L_0\cdot
v_\alpha=\Dslzr(L_0)+\sum_{\alpha\in\Delta_s}\alpha\,v_\alpha$ is an infinite sum, which is not in $\VVsl$, contradicting the fact that $\Dsl$ is a derivation from $\LLsl$ to $\VVsl$. Then Claim \ref{finsumclaim} follows.

\begin{clai}\label{H1LQsleq0}
$H^1(\LLsl,\QQsl)=0$.
\end{clai}

The exact sequence $0\rar\frak{h}^s\rar\LLsl\rar\LLsl/\frak{h}^s\rar0$ induces an exact sequence of low degree in the Hochschild-Serre spectral sequence
\begin{eqnarray*}
0\lrar H^1(\LLsl/\frak{h}^s,{\QQsl}^{\frak{h}^s})\lrar H^1(\LLsl,\QQsl)\lrar {H^1(\frak{h}^s,\QQsl)}^{\LLsl/\frak{h}^s}.
\end{eqnarray*}
Lemma \ref{WittVirLiebialg} forces $H^1(\LLsl/\frak{h}^s,{\QQsl}^{\frak{h}^s})=0$. Since $\LLsl/\frak{h}^s\cong\frak{w}$ and $\QQsl\cong\frak{w}\otimes\frak{w}$, ${H^1(\frak{h}^s,\QQsl)}^{\LLsl/\frak{h}^s}$ can be embedded into ${\rm Hom}_{U(\frak{w})}(\frak{h}^s,\frak{w}\otimes\frak{w})$, which can be easily proved to be zero. Then Claim \ref{H1LQsleq0} follows.

\begin{clai}\label{H1LHLCsleq0}
$H^1(\LLsl,\HHsl/\LCsl)=0$.
\end{clai}

Denote the subalgebra spanned by $\{L_n,\,M_n,\,\c\,|\,n\in\Z\}$ of $\LLsl$ as $\WWsl$, $\JJsl=\WWsl\otimes\frak{h}^s+\frak{h}^s\otimes\WWsl$, $\WCsl=\WWsl\otimes\frak{C}^s_{\lambda}+\frak{C}^s_{\lambda}\otimes\WWsl$ and $Y_{s}=\{Y_{n+s}\,|\,n\in\Z\}$. This claim of the case $s\!=\!\frac12$ follows immediately from the following Subclaim \ref{H1WJWCsl=0} and Subclaim \ref{H1WYYs=0}.

\begin{Subclai}\label{H1WJWCsl=0}
$H^1(\WWsl,\JJsl/\WCsl)=0$.
\end{Subclai}

This subclaim can be found in Remark 1 of \cite{LPZ-JA12}.

\begin{Subclai}\label{H1WYYs=0}
$H^1(\WWsl,Y_{s}\otimes Y_{s})=0$ for the case $s=\frac12$.
\end{Subclai}

This subclaim can be proved similar to the proof of Theorem 4.5 (i) of  \cite{LPZ-JA12}.
According to Subclaim \ref{H1WJWCsl=0} and Subclaim \ref{H1WYYs=0}, we know that
$\vpsl(L_{n})=\vpsl(M_{n})=0$ for any $n\in\Z$ and $\vpsl\in H^1(\LLsl,\HHsl/\LCsl)$ when $s=\frac12$, which implies $\vpsl(Y_{\frac12})=0$. Then Claim \ref{H1LHLCsleq0} holds for the case $s=\frac12$.

Similar to the proof of Claim 3 of \cite{FLL}, we can prove this claim also holds for the case $s=0$.

\begin{clai}\label{H1LLCslCong}
\begin{eqnarray*}
\begin{array}{lll}
&&H^1(\LLsl,\LCsl)\cong
\left\{\begin{array}{lllll}
\frak{D}^\natural_1&\mbox{\rm{if}\ \ }s=\frac12,\ \,\lambda=-1,\\[3pt]
\frak{D}^\natural_2&\mbox{\rm{if}\ \ }s=\frac12,\ \,\lambda=-2,\\[3pt]
\frak{D}^\natural_3&\mbox{\rm{if}\ \ }s=\frac12,\ \,\lambda\notin\mathcal{S}_\frac12,\\[3pt]
\frak{d}^\natural_1&\mbox{\rm{if}\ \ }s=0,\ \,\lambda=0,\\[3pt]
\frak{d}^\natural_2&\mbox{\rm{if}\ \ }s=0,\ \,\lambda=-1,\\[3pt]
\frak{d}^\natural_3&\mbox{\rm{if}\ \ }s=0,\ \,\lambda=-2,\\[3pt]
\frak{d}^\natural_4&\mbox{\rm{if}\ \ }s=0,\ \,\lambda=1,\\[3pt]
\frak{d}^\natural_5&\mbox{\rm{if}\ \ }s=0,\ \,\lambda\notin\mathcal{S}_0.
\end{array}\right.
\end{array}
\end{eqnarray*}
\end{clai}

This claim follows from Lemmas \ref{12Zcenterlamm}, \ref{OrderTensor}, \ref{TwderTensor} and Claims \ref{H1LQsleq0}, \ref{H1LHLCsleq0}.\vskip6pt

By now we have proved Proposition \ref{propCohLlsVls}.\QED\vskip6pt

The following lemma is still true for $\LLsl$ by employing the technique of Lemma 2.5 in \cite{HLS}.

\begin{lemm}\label{xv1mtauv1mtau}
Suppose $v\in\VVsl$ such that $x\cdot v\in {\rm Im}(1-\tau)$ for all
$x\in\LLsl.$ Then there exists some $u\in {\rm Im}(1-\tau)$ such that $v-u\in\frak{C}_s\otimes\frak{C}_s$.
\end{lemm}

\ni{\bf Proof of Theorem \ref{maintheo}}\ \ \,This theorem follows from Lemmas \ref{OrderTensor}, \ref{TwderTensor}, \ref{xv1mtauv1mtau} and Proposition \ref{propCohLlsVls}.\QED


\begin{thebibliography}{9999}\parskip0pt\lineskip4pt\small

\bibitem{D1} V.G. Drinfeld, Constant quasiclassical solutions of the
Yang-Baxter quantum equation, {\it Soviet Math. Dokl.} {\bf28}(3)
(1983), 667--671.

\bibitem{D2} V.G. Drinfeld, Quantum groups, in: {\it Proceeding of the
International Congress of Mathematicians}, Vol.~1, 2, Berkeley,
Calif.~1986, Amer.~Math.~Soc., Providence, RI, 1987, 798--820.

\bibitem{FLL} H. Fa, Y. Li, J. Li, The Schr\"{o}dinger-Virasoro type Lie bialgebra: a
twisted case, {\it Front. Math. China} {\bf6} (4) (2011) 641--657.

\bibitem{FLX} H. Fa, J. Li, B. Xin, Lie superbialgebra structures on the centerless
twisted $N\!=\!2$ superconformal algebra, {\it Algebra Colloq.} {\bf18} (3) (2011) 361--372.

\bibitem{GJP-AC} S. Gao, C. Jiang, Y. Pei, Structure of the extended
Schr\"{o}dinger-Virasoro Lie algebra, {\it Alg. Colloq.} {\bf 16} (2009), 549--566.

\bibitem{HLS} J. Han, J. Li, Y. Su, Lie bialgebra structures on the
Schr\"{o}dinger-Virasoro Lie algebra, {\it J. Math. Phys.} {\bf 50} (2009), 083504.

\bibitem{H1} M. Henkel, Schr\"{o}dinger invariance and strongly
anisotropic critical systems, {\it J. Stat. Phys.} {\bf 75} (1994),
1023--1029.

\bibitem{H2} M. Henkel, Phenomenology of local scale invariance: from
conformal invariance to dynamical scaling, {\it Nucl. Phys. B} {\bf
641} (2002), 405--410.

\bibitem{HU} M. Henkel, J. Unterberger, Schr\"{o}dinger invariance and
space-time symmetries, {\it Nucl. Phys. B} {\bf 660} (2003), 407--412.

\bibitem{JW-arXiv12} Q. Jiang, S. Wang, Derivations and automorphism groups of
the original deformative Schr\"{o}dinger-Virasoro algebras, arXiv:1209.3164v1.

\bibitem{LJ} J. Li, 2-cocycles of twisted deformative Schr\"{o}dinger-Virasoro algebra,
{\it Comm. Algebra} {\bf 40} (2012), 1933--1950.

\bibitem{LCZ} D. Liu, L. Chen, L. Zhu, Lie superbialgebra structures on the $N\!=\!2$
superconformal Neveu-Schwarz algebra, {\it J. Geometry and Phys.} {\bf 62} (2012), 3826--831.

\bibitem{LPZ-JA12} D. Liu, Y. Pei, L. Zhu, Lie bialgebra structures on the twisted
Heisenberg-Virasoro algebra, {\it J. Alg.} {\bf 359} (2012) 35--48.

\bibitem{LSZ} J. Li, Y. Su, L. Zhu, 2-cocycles of original deformative
Schr\"{o}dinger-Virasoro algebras, {\it Science in China Series A} {\bf 51} (2008), 1989--1999.

\bibitem{LS-JMP} J. Li, Y. Su, Representations of the Schr\"{o}dinger-Virasoro
algebras, {\it J. Math. Phys.} {\bf 49} (2008), 053512.

\bibitem{LS-arXiv08} J. Li, Y. Su, The derivation algebra and automorphism group
of the twisted Schr\"{o}dinger-Virasoro algebra, arXiv:0801.2207v1,
(2008).

\bibitem{LSX} J. Li, Y. Su, B, Xin, Lie bialgebras of a family of Block
type, {\it Chinese Annals of Math. (Series B)} {\bf 29} (2008), 487--500.

\bibitem{M1} W. Michaelis, A class of infinite-dimensional Lie bialgebras
containing the Virasoro algebras, {\it Adv. Math.} {\bf107} (1994),
365--392.

\bibitem{M2} W. Michaelis, Lie coalgebras, {\it Adv. Math.} {\bf38}
(1980), 1--54.

\bibitem{M3} W. Michaelis, The dual Poincare-Birkhoff-Witt theorem, {\it
Adv.Math.} {\bf57} (1985), 93--162.

\bibitem{NT} S.H. Ng, E.J.~Taft, Classification of the Lie bialgebra
structures on the Witt and Virasoro algebras, {\it J. Pure Appl. Alg.} {\bf151} (2000), 67--88.

\bibitem{RU} C. Roger, J. Unterberger, The Schr\"{o}dinger-Virasoro Lie
group and algebra: representation theory and cohomological study, {\it Ann. Henri Poincar$\acute{e}$} {\bf 7} (2006),  1477--1529.

\bibitem{SS} G. Song, Y. Su, Lie bialgebras of generalized Witt type,
{\it Science in China: Series A} {\bf49} (2006), 533--544.

\bibitem{T} E.J. Taft, Witt and Virasoro algebras as Lie bialgebras, {\it
J. Pure Appl. Alg.} {\bf87} (1993), 301--312.

\bibitem{TZ-JA09}  S. Tan, X. Zhang, Automorphisms and Verma modules for generalized
Schr\"{o}dinger-Virasoro algebras, {\it J. Algebra} {\bf 322} (2009), 1379--1394.

\bibitem{U} J. Unterberger, On vertex algebra representations of the
Schr\"{o}dinger-Virasoro Lie algebra, {\it Nuclear Physics B} {\bf 823} (2009), 320--371.

\bibitem{WLXin} W. Wang, J. Li, B. Xin, Central extensions and derivations of
generalized Schr\"{o}dinger-Virasoro algebras, {\it Alg. Colloq.} {\bf 19}  (2012), 735--744.

\bibitem{WLXu} W. Wang, J. Li, Y. Xu, Derivations and automorphisms of twisted
deformative Schr\"{o}dinger-Virasoro Lie algebras, {\it Comm. Algebra} {\bf 40} (2012), 3365--3388.

\bibitem{WSS06} Y. Wu, G. Song, Y. Su, Lie bialgebras of generalized
Virasoro-like type, {\it Acta Mathematica Sinica, English Series} {\bf22} (2006), 1915--1922.

\bibitem{WSS07} Y. Wu, G. Song, Y. Su, Lie bialgebras of generalized Witt
type. II, {\it Comm. Algebra} {\bf35} (6) (2007), 1992--2007.

\bibitem{XL} Y. Xu, J. Li, Lie bialgebra structures on the extended affine Lie algebra
$\widetilde{\frak{sl}_2(\C_q)}$, {\it J. Pure and Applied Alg.} {\bf217} (2013), 364--376

\bibitem{XSS} B. Xin, G. Song, Y. Su, Hamiltonian type Lie bialgebras,
{\it Science in China A} {\bf50} (2007), 1267--1279.

\bibitem{Y-HS09} H. Yang, Y. Su, Lie bialgebras structures on the Ramond $N\!=\!2$
super-Virasoro algebras, {\it Chaos, Solitons and Fractals} {\bf 40} (2009),
661--671.

\bibitem{Y-XS08} X. Yue, Y. Su, Lie bialgebra structures on Lie algebras of
generalized Weyl type, {\it Comm. Algebra} {\bf36}(4)(2008), 1537--1549.

\bibitem{ZTL} X. Zhang, S. Tan, H. Lian, Whittaker modules for the Schr\"{o}dinger-Witt
algebra, {\it J. Math. Phys.} {\bf 51} (2010), 083524.







\end{thebibliography}
\end{document}